\documentclass[11pt]{article}
\usepackage{amssymb}
\usepackage{mathrsfs}
\usepackage{latexsym,amsmath,amssymb,amsfonts,epsfig,graphicx,cite,psfrag}
\usepackage{eepic,color,colordvi,amscd}

\newcommand{\qed}{\hfill $\Box $}
\newcommand{\pf}{\noindent {\bf Proof.} }
\newtheorem{theorem}{Theorem}[section]
\newtheorem{lemma}[theorem]{Lemma}

\newtheorem{conjecture}[theorem]{Conjecture}

\begin{document}

\title{Matching stability for 3-partite 3-uniform hypergraphs\thanks{Supported by  National Key Research and Development Program of China 2023YFA1010203 and the National Natural
Science Foundation of China under grant No.12271425}}

\author{Hongliang Lu and Xinxin Ma\\School of Mathematics and Statistics\\
Xi'an Jiaotong University\\
Xi'an, Shaanxi 710049, China\\
\smallskip\\
}
\date{}

\maketitle

\date{}

\maketitle

\begin{abstract}
Let $n,k,s$ be three  integers such that $k\geq 2$ and $n\geq s\geq 1$. Let $H$ be a $k$-partite $k$-uniform hypergraph with $n$ vertices in each class.
Aharoni (2017) showed that if $e(H)>(s-1)n^{k-1}$, then $H$ has a matching of size $s$.
 In this paper, we gve a stability result for 3-partite 3-uniform hypergraphs: if  $G$ is a $3$-partite $3$-uniform hypergraph  with $n\geq 162$ vertices in each class, $e(G)\geq (s-1)n^2+3n-s$ and $G$ contains no matching of size $s+1$, then $G$ has a vertex cover of size $s$. Our bound is also tight.
 \end{abstract}

\section{Introduction}
Let $n,k$ be two positive integers and let $[n]:=\{1,\ldots,n\}$. For a set
$S$, let ${S\choose k}:=\{T\subseteq S: |T|=k\}$. A {\it hypergraph} $H$
consists of a vertex set $V(H)$ and an edge set $E(H)$ whose members
are subsets of $V(H)$. Write $e(H):=|E(H)|$.  A hypergraph $H$ is {\it $k$-uniform} if $E(H)\subseteq {V(H)\choose k}$, and a $k$-uniform hypergraph is also
called a {\it $k$-graph}. Thus, 2-graphs are graphs.
A $k$-graph is $k$-$partite$ if there exists a
partition of the vertex set $V(H)$ into sets $V_1, \cdots, V_k$ (called {\it partition  classes}) such that for any $f\in E(H)$, $|f\cap V_i|=1$
for $i\in [k]$.

Let $H$ be a $k$-graph and $T\subseteq V(H)$. The neighbour of $T$ in $H$ denoted by $N_H(T)$ is $\{S\subseteq V(H)\ |\ S\cup T\in E(H)\}$. When there is no confusion, we also use $N_H(T)$ to denote the hypergraph with vertex set $V(H)\setminus T$ and edge set $N_H(T)$.   The {\it degree}
of $T$ in $H$, denoted by $d_H(T)$, is the number of edges of $H$
containing $T$, i.e., $d_H(T):=|N_H(T)|$.
 Let $l$ be a nonnegative integer; then $\delta_l(H):=
\min\{d_H(T): T\in {V(H)\choose l}\}$ denotes the {\it minimum
$l$-degree} of $H$. $\delta_0(H)$ is the number of edges in $H$, $\delta_1(H)$ is often called the  {\it minimum
vertex} degree of $H$.
A {\it matching} in $H$ is a set of pairwise disjoint
edges of $H$, and it is \emph{perfect} if the union of all edges in
the matching is $V(H)$. A \emph{maximum matching} is a matching of the maximum number of edges.
The size of a maximum matching in  $H$ is called \emph{matching number} of $H$, denoted by $\nu(H)$. 
 A vertex subset $S\subseteq V(H)$ is called a \emph{vertex cover} of $H$   if every edge in $H$ intersects $S$. A \emph{minimum vertex cover} of $H$ is a vertex cover having the smallest possible number of vertices.
  We use $\tau(H)$ to
denote the size of a minimum vertex cover of  $H$. Let $\mathcal{F}= \{F_1, ... ,F_m\}$ be a family of hypergraphs on the same vertex set.
A set of $m$ pairwise disjoint edges is called a \emph{rainbow matching} for $\mathcal{F}$ if each edge is from a different $F_i$.
If such edges exist, then we also say that $\mathcal{F}$ admits a rainbow matching.

Erd\H{o}s and Gallai \cite{EG59} determined $e(G)$ for a graph $G$ to contain a matching of given size, and  Erd\H{o}s~\cite{Er65} conjectured
the following generalization to $k$-graphs for $k\ge 3$: The threshold on $\delta_0(H)$ for a $k$-graph $H$ on $n$ vertices to contain a matching of size $m$ is
$$\max\left\{{km-1\choose k},{n\choose k}-{n-m+1\choose k}\right\}+1.$$
For recent related results on this conjecture, we refer the readers to \cite{AF12,Fr13,Fr17,Fr17-2,FK19,FK18,FRR12,HLS,KOT14,LM14}.

There has also been more interest in matchings in $k$-partite $k$-graphs.
For a $k$-partite $k$-graph $H$, a set $T\subseteq V(H)$ is said to be {\it
legal} if $|T\cap V_i|\leq 1$ for all $i\in [k]$, and {\it balanced} if $|T\cap V_i|=|T\cap V_j|$ for all $i,j\in [k]$. Thus, if $T$ is
not legal in $H$ then $d_H(T)=0$. For integers $0\le l\le k-1$, let $\delta_l(H):= \min\{d_H(T):
T\in {V(H)\choose l} \mbox{ and $T$ is legal}\}$. Aharoni, Georgakopoulos
and Spr\"{u}ssel \cite{AGS09} proved a result about codegrees for perfect matchings in
$k$-partite $k$-graphs, answering a question of K\"{u}hn and Osthus \cite{KO06}. Berge \cite{Berge75} proved that every $k$-partite $k$-graph $H$ may be decomposed into $\Delta(H)$ edge-disjoint matchings, where $\Delta(H)$  denotes the maximum vertex degree.
\begin{theorem} [Berge, \cite{Berge75}] \label{Berge}
Let $n,k,m$ be three positive integers such that $n\geq m$, and let $H$ be a complete $k$-partite $k$-graph with $n$ vertices  in each partition class.
Then $E(H)$ can be decomposed into $n^{k-1}$ edge-disjoint perfect matchings.
\end{theorem}

In term of Theorem \ref{Berge}, Aharoni and Howard  \cite{AH17}   obtained the following result.
\begin{theorem} [Aharoni and Howard, \cite{AH17}] \label{Aharoni17}
Let $n,k,m$ be three positive integer such that $n\geq m$, and let $H$ be a $k$-partite $k$-graph in each partition class.
If $e(H)> (m-1)n^{k-1}$, then $\nu(H)\geq m$.
\end{theorem}


In this paper, we study the stability version for matchings in 3-partite 3-graphs.
\begin{theorem}\label{main-theorem}
Let $n,m$ be two integers such that $n>\max\{m,162\}$. If $H$ is a
$3$-partite $3$-graph with $n$ vertices in each partition class,
$e(H)\geq (m-1)n^2+3n-m$ and $\nu(H)=m$, then $\tau(H)= m$.
\end{theorem}

The bound in   Theorem \ref{main-theorem} is tight. Let $n,k,m$ be three integers such that $n\geq m+1\geq 2$ and $k\geq 3$. Let $V_1,\ldots, V_k$ be $k$ pair-disjoint sets of order $n$ and let $W_1\subseteq V_1$ such that $|W_1|= m-1$. We choose $k$ vertices $v_1,\ldots,v_k$ such that $v_1\in V_1\backslash W_1$ and   $v_i\in V_i$ for $2\leq i\leq k$. Consider the $k$-partite $k$-graph
$H_k(n,m)$  with partition classes $V_1,\ldots,V_k$ and edge set $E(H_k(n;m))$ consisting of all legal $k$-subsets of $V(H_k(n;m))$ intersecting $W_1$ or containing $\{v_2,\ldots,v_k\}$ or containing $v_1$ and intersecting $\{v_2,\ldots,v_k\}$.
Clearly,
$\nu(H_k(n,m))=m$ (when $n\ge m$), $\tau(H_k(n,m))=m+1$ and
\[
e(H_k(n,m))=(m-1)n^{k-1}+(n^{k-1}-(n-1)^{k-1})+n-m.
\]

Based the construction of $H_k(n,m)$, we would like to propose the following conjecture.
\begin{conjecture}
Let $n,m,k$ be three positive integers such that $n$ is sufficiently large, $n>m$ and $k\geq 4$. If $H$ is a
$k$-partite $k$-graph with $n$ vertices in each partition class and
$e(H)\geq e(H_k(n,m))+1$ and $\nu(H)=m$, then $\tau(H)= m$.
\end{conjecture}


Let $H$ be a hypergraph, $A,B$ disjoint subsets of $V(H)$, and $x\in
V(H)\setminus(A\cup B)$. We use $L_x(A,B)$ to denote the {\it link
graph} of $x$ with respect to $A,B$, which is a bipartite graph
with partition classes $A$ and $B$ such that, for $u\in A$ and $v\in
B$, $uv$ is an edge in $L_x(A,B)$ if and only if
$\{u,v,x\}\in E(H)$. 
 For any $S\subseteq V(H)\setminus(A\cup B)$, we let $L_S(A,B):=\bigcup_{x\in S}L_x(A,B)$ and $L_S(A):=L_S(A,A)$. We use $e_H(A,B)$ to denote the number of edges intersecting $A$ and $B$.


We end this introductory section with the following convention.
Let $H$ be a hypergraph. 
For any $M\subseteq E(H)$,
let $V(M):=\bigcup_{f\in M}f$. For any $T\subseteq V(H)$, we use $H[T]$ to denote the
hypergraph with vertex set $T$ and edge set $\{f\in E(H): f\subseteq T\}$, and let $H-T:=H[V(H)\setminus T]$.

\section{Preliminaries}

In this section, we present a few lemmas for our proof of Theorem~\ref{main-theorem}.
We begin with several results  from \cite{LM14}, first of which will serve as the base of our proof.

Lo and Markstr\"om  gave the following lemma (Lemma 5.5 in \cite{LM14}).
\begin{lemma} [Lo and Markstr\"om, \cite{LM14}] \label{extension}
Let $H$ be a 3-partite 3-graph, $M$ a matching in $H$, and $S$ a legal 3-subset of $X=V(H)\backslash V(M)$. For any $\{f_1,f_2\}\in {M\choose 2}$, if $e(L_S(f_1,f_2))\geq 5$, then $L_S(f_1,f_2)$ has a perfect matching.
\end{lemma}

When $M$ is a maximum matching, we can  apply Lemma~\ref{extension} to bound the number of pairwise disjoint legal
3-subsets $S$ of $X$ for which $L_S(f,g)$ has a large number of edges for any $\{f,g\}\in {M\choose 2}$.
%
In our proof, we also need the following result.
\begin{theorem} [Aharoni and Howard, \cite{AH17}] \label{Aharoni17-2}
Let $n,m$ be two positive integers such that $n\geq m$, and let $G_1,\ldots,G_m$ be $m$ bipartite graphs with the same vertex set $A\cup B$, where $|A|=|B|=n$.
If $e(G_i)> (m-1)n$, then $\{G_1,\ldots,G_n\}$ admits a rainbow matching.
\end{theorem}


The proof in \cite{EKR} uses an operation called 'shifting', which has since become a main tool
in the area. It is an operation on a hypergraph $G$, defined with respect to a specific linear
ordering '$<$' on its vertices. For $x<y$ in $V(G)$, we define $S_{xy}(e)=(e\cup \{x\}) \setminus \{y\}$ if $e\in E(G)$. $x\in e$ and $y\notin e$,
provided $(e\cup \{x\})\setminus \{y\} \notin E(G)$; otherwise let $S_{xy}(e)=e$. We also write $S_{xy}(G) =\{S_{xy}(e)\ |\ e\in E(G)\}$.
Given a $k$-partite $k$-graph $G$ and a linear order on its vertex class, a $k$-partite
shifting is a shifting $S_{xy}$ where $x$ and $y$ belong to the same vertex class. $G$ is said to be
\emph{partitely shifted} if $S_{xy}(G)=G$ for all pairs $x<y$ that belong to the same vertex class.

Erd\H{o}s and Gallai \cite{EG61} showed that the shifting operation does not increase the matching number of a
hypergraph.
\begin{lemma}
Let $H$ be a $k$-graph. For any two different vertices $x,y$ of $H$, $\nu(S_{xy}(H))\leq \nu(H)$.
\end{lemma}

Now we give the stability result for intersecting family.
\begin{lemma}\label{3-inter-family}
Let $n_1, n_2,n_3$ be three integers such that $\min \{n_1, n_2, n_3\}\geq 5$. Let $H$ be $3$-partite 3-graph with vertex partition $(V_1,V_2,V_3)$, where $|V_i|=n_i$ for $1\leq i\leq 3$. If $e(H)\geq n_1+n_2+n_3-1$ and $\nu(H)=1$, then $\tau(H)= 1$.
\end{lemma}
\pf By contradiction. Suppose that $\tau(H)\geq 2$. Then $E(H-v)\neq \emptyset$ for any $v\in V(H)$.  Let $f_0\in E(H)$. Since $e(H)\geq n_1+n_2+n_3-1$, and every edge of $H$ intersects $f_0$, then $\sum_{x\in f_0} d_{H}(x)\geq n_1+n_2+n_3-1$. Thus $f_0$ contains a vertex, saying $u_1\in V_1$ such that $d_{H}(u_1)\geq 5$. Let $f_1=\{v_1,v_2,v_3\}\in E(G-u_1)$ such that $v_i\in V_i$ for all $1\leq i\leq 3$.
Since $d_H(u_1)\geq 5$ and $\nu(H)=1$, by symmetry, we may assume there exist two vertices $a_3,b_3$ in $V_3$ such that $\{u_1,v_2,a_3\},\{u_1,v_2,b_3\}\in E(H)$.  Write $f_2:=\{u_1,v_2,a_3\}$ and $f_3:=\{u_1,v_2,b_3\}$.

\medskip
\textbf{Claim 1.}~For any $f\in E(H)$, $|f\cap \{u_1,v_2,v_3\}|\geq 2$.
\medskip

Firstly, consider $x\in V_2\setminus\{v_2\}$. We claim $N_{H}(x)\subseteq \{\{u_1,v_3\}\}$. Otherwise, there exists  $g\in N_{H}(x)\setminus\{\{u_1,v_3\}\}$. If $u_1\notin g$, then either
$\{f_2,g\cup \{x\}\}$ or $\{f_3,g\cup \{x\}\}$ is a matching of size two in $H$, contradicting to that $\nu(H)=1$; else if $u_1\in g$ and $v_3\notin g$, then  $\{g\cup \{x\},f_1\}$ is a matching of size two in $H$, a contradiction. Since $\tau(H)\geq 2$, there exists an edge  $f_4=\{u_1,a_2,v_3\}\in E(H)$ such that $a_2\in V_2-v_2$.

Secondly, consider $y\in V_1\setminus\{u_1\}$. We claim $N_{H}(y)\subseteq \{\{v_2,v_3\}\}$. Otherwise,   there exists $g\in N_{H}(y)\setminus \{\{v_2,v_3\}\}$. Since $\nu(H)=1$, $(g\cup \{y\})\cap f_1\neq \emptyset$.  If $v_2\notin g$, then either $\{g\cup \{y\}, f_2\}$ or $\{g\cup \{y\}, f_3\}$ is a matching of size two in $H$, a contradiction; else if $v_2\in g$ and $v_3\notin E(H)$, then   $\{g\cup \{y\},f_4\}$  is a matching of size two in $H$,  a contradiction.

Finally, consider $z\in V_3-v_3$. We claim $N_{H}(z)\subseteq \{\{u_1,v_2\}\}$. Otherwise,   there exists  $g\in N_{H}(z)\setminus\{\{u_1,v_2\}\}$. Since $\nu(H)=1$, then $g\cap \{u_1,v_2\}\neq \emptyset$, otherwise, either $\{g\cup\{z\}, f_2\}$ or $\{g\cup\{z\}, f_3\}$ is a matching of size two in $H$, a contradiction. If $u_1\notin g$ and $v_2\in g$, then $\{f_4,g\cup \{z\}\}$ is a matching of size two in $G$, a contradiction; else if $u_1\in g$ and $v_2\notin g$, then $\{f_1,g\cup \{z\}\}$ is a matching of size two in $H$, a contradiction. This completes the proof of Claim 1.

By Claim 1, every edge of $G$ intersects at least two vertices of $\{u_1,v_2,v_3\}$. So we have $e(H)\leq n_1+n_2+n_3-2$, a contradiction. This completes the proof. \qed

\begin{lemma}\label{matching-2size}
Let $n\geq 2$ be an integer.
Let $G_1$ and $G_2$ be two bipartite graphs with $n$ vertices in each class. If for $i\in [2]$, $e(G_i)>0$ and $e(G_1)+e(G_2)>2n$, then $\{G_1,G_2\}$ admits a rainbow matching.
\end{lemma}

\pf Suppose that $e(G_1)\leq e(G_2)$. Then $e(G_2)>n$. Let $uv\in E(G_1)$. If $G_2-\{u,v\}$ contains an edge $f$, then $\{uv,f\}$ is a rainbow matching of $\{G_1,G_2\}$. So we may assume that $\{u,v\}$ is a vertex cover of $G_2$. Thus we may infer that $d_{G_2}(u)\geq 2$ and $d_{G_2}(v)\geq 2$. Consider $e(G_1)\geq 2$. Then we may choose an edge $f'\in E(G_1)-uv$. Without loss of generality, suppose that $v\notin f'$. Since $d_{G_2}(v)\geq 2$ and $G_2$ is a bipartite graph, then there exists an edge $f''\in E(G_2)$ such that $v\in f''$ and $f'\cap f''=\emptyset$. It follows that $\{f',f''\}$ is a rainbow matching of $\{G_1,G_2\}$. Next we may assume that $E(G_1)=\{uv\}$, which implies that $e(G_2)\geq 2n$, which contradicts to that $\{u,v\}$ is a vertex cover of $G_2$. This completes the proof.

Next we give the stability lemma for smaller matching.

\begin{lemma}\label{matching-small}
Let $n,m$ be two integers such that $m\geq 2$ and $n\geq 12m$. Let $H$ be a 3-partite 3-graph with $n$ vertices in each class. If $e(H)\geq (m-1)n^2+3n-m$ and $\nu(H)=m$, then $\tau(H)=m$.

\end{lemma}
\pf We prove the following claim.

\medskip
\textbf{Claim 1.}~Let $s< n/2$ be an integer and let $G$ be a $3$-partite $3$-graph with $n$ vertices in each class.  If $G$ contains a vertex $v$ such that $\nu(G-v)=s$ and $d_{G}(v)>2sn$, then $\nu(G)\geq s+1$.  \medskip

Let $M$ be a matching of size $s$ in $G-v$. Note that the number of edges containing $v$ and intersecting $V(M)$ is at most $2sn$. So there exists an edge $f$ such that $v\in f$ and $f\cap V(M)=\emptyset$. Then $M\cup \{f\}$ is a matching of size $s+1$ in $G$. This completes the proof of Claim 1.

\medskip
\textbf{Claim 2.}~$H$ contains $m-1$ vertices, say $v_1,\ldots, v_{m-1}$ such that $d_{H-A_i}(v_i)> 2(m-i)n$ and $\nu(H-A_{i+1})=m-i+1$, where $A_1:=\emptyset$ and $A_i:=\{v_1,\ldots,v_{i-1}\}$ for $2\leq i\leq m$.
\medskip

Let $M_0$ be a maximum  matching of $H$. Then we have $|M_0|=m$. Note that
\[
\sum_{v\in V(M_0)}d_{H}(v)\geq e(H)> (m-1)n^2+2n,
\]
Hence there exists $v_1\in V(M_0)$ such that
\[
 d_H(v_1)> \frac{1}{3m}((m-1)n^2+2n)> 2mn,
\]
where the second inequality holds since $n\geq 12m$ and $m\geq 2$.
Since $\nu(H)=m$, then by Claim 1, we have $\nu(H-v_1)=m-1$.  Recall that $e(H)\geq (m-1)n^2+3n-m$. So we have $e(H-v_1)\geq  (m-2)n^2+3n-m$. Suppose that we have found $v_1,\ldots,v_r$ such that for $1\leq i\leq r\leq m-2$, the following statements hold.
\begin{itemize}
  \item  [(i)] $\nu(H-\{v_1,\ldots,v_{i}\})=m-i$;
  \item [(ii)] $d_{H-A_i}(v_i)> 2(m-i)n$ for $1\leq i\leq r$, where $A_1=\emptyset$ and $A_i=\{v_1,\ldots,v_{i-1}\}$ for $2\leq i\leq r$.
\end{itemize}
Suppose that $r\leq m-2$. Note that $e(H-\{v_1,\ldots,v_{r}\})>(m-1-r)n^2+2n$. Let $M_r$ be a matching of size $m-r$ in $H-\{v_1,\ldots,v_{r}\}$. Then we have
\[
\sum_{v\in V(M_r)}d_{H}(v)> (m-1-r)n^2+2n.
\]
Hence there exists $v_{r+1}\in V(M_r)$ such that
\[
 d_H(v_{r+1})> \frac{1}{3(m-r)}((m-1-r)n^2+2n)> 2(m-r)n,
\]
where the second inequality holds since $n\geq 12m$ and $r\leq m-2$.
Moreover, by Claim 1, we have $\nu(H-\{v_1,\ldots,v_{r+1}\})=m-1-r$. Continuing the process, we may find the desired $m-1$ vertices. This completes the proof of Claim 2.

By Claims 1 and 2, we may choose  $v_1,\ldots,v_{m-1}$ such that for $1\leq j\leq m-1$, $\nu(H-\{v_1,\ldots,v_{j}\})=m-j$.  Write $A:=\{v_1,\ldots,v_{m-1}\}$ and $n_i:=n-|V_i\cap \{A\}|$ for $1\leq i\leq 3$. Note that
$\min \{n_1,n_2,n_3\}\geq n-m\geq 5$.
One can see that $\nu(H-A)=1$ and 
\[
e(H-A)\geq 3n-m= n_1+n_2+n_3-1.
\]
By Lemma \ref{3-inter-family}, $H-A$ has a vertex cover  of size one denoted by $A'$.  Then $A\cup A'$ is a vertex cover of size $m$ in $H$. This completes the proof. \qed

\section{Matching Stability for Partitely Shifted  Graphs}

In Sections 3 and 4,  let $V_1,V_2,V_3$ be three pair-disjoint sets such that $V_1:=\{a_1,\ldots,,a_n\}$, $V_2:=\{b_1,b_2\ldots,b_n\}$ and $V_3:=\{c_1,\ldots,c_n\}$, where $a_1<a_2\ldots<a_n$, $b_1<b_2\ldots<b_n$ and $c_1<c_2\ldots<c_n$.
\begin{lemma}\label{Stab-lem}
Let $n,m$ be two integers such that $ n\geq m+2$ and $m\geq 2$. Let $G$ be partitely shifted  3-partite 3-graph with vertex partition $(V_1,V_2,V_3)$.  If $e(G)\geq (m-1)n^2+3n-m$ and $\nu(G)=m$, then $\tau(G)=m$.
\end{lemma}

\pf By Lemma \ref{matching-small}, we may assume that $n\leq 12m$. By induction on $m$, suppose that the result holds for smaller $m$.   Without loss of generality, we may assume that $G$ is a maximal graph with $\nu(H)=m$. Since $G$ is  partitely shifted, then we have $N_{G}(a_n)\subseteq \cdots\subseteq N_G(a_1)$, $N_{G}(b_n)\subseteq \cdots\subseteq N_G(b_1)$ and $N_{G}(c_n)\subseteq \cdots\subseteq N_G(c_1)$.

\medskip
\textbf{Claim 1.}~$N_{G}(a_n)= \cdots=N_G(a_{m+1})$, $N_{G}(b_n)=\cdots= N_G(b_{m+1})$ and $N_{G}(c_n)= \cdots=N_G(c_{m+1})$.
\medskip

By symmetry, it is sufficient for us to show that $N_{G}(a_{i})=N_{G}(a_{i+1})$ for $m+1\leq i\leq n-1$. By contradiction, suppose that there exists $i\geq m+1$ such that  $N_{G}(a_{i})\neq N_{G}(a_{i+1})$. Let $g\in N_{G}(a_{i})\backslash N_{G}(a_{i+1})$. Since $G$ is maximal, $E(G)\cup \{g\cup \{a_{i+1}\}\}$ has a matching $M$ of size $m+1$. Since $\nu(G)=m$, we have $g\cup \{a_{i+1}\}\in M$. So there exists $j\leq m+1$ such that $a_j\notin V(M)$. Recall that $N_{G}(a_i)\subseteq N_{G}(a_j)$. Thus $(M\backslash \{g\cup \{a_{i+1}\}\})\cup \{g\cup \{a_j\}\}$ is a matching of size $m+1$ in $G$, a contradiction. This completes the proof of Claim 1.

If there exists $v\in \{a_{m+1},b_{m+1},c_{m+1}\}$ such that $d_{G}(v)=0$, then by stability of $G$, one can see that $\tau(G)\leq m$. So by Claim 1, we may assume that $G$ contains no isolated vertices.

 \medskip
\textbf{Claim 2.}~For any $v\in V(G)$, $\nu(G-v)=m$.
\medskip

Suppose that $G$ contains a vertex, saying $a_1$ such that $\nu(G-a_1)=m-1$. Consider $N_{G}(a_1)= N_G(a_{m+1})$. Since $\nu(G)=m$, then we have $\nu(N_G(a_1))=m$. Note that $N_G(a_1)$ is a bipartite graph. By K\"onig's Theorem, $N_{G}(a_1)$ has a vertex cover of size $m$, which is also a vertex cover of $G$. Hence we may suppose that $N_{G}(a_1)\backslash N_G(a_{m+1})\neq \emptyset$. Let $G'$ be spanning subgraph of $G$ with edge set $E(G)\backslash \{g\cup {a_1}\ |\ g\in N_{G}(a_1)\backslash N_G(a_{m+1})\}$. One can see that $e(G')\geq (m-2)n^2+3n-m+1$ since $\delta_1(G)\geq 1$. Now by induction hypothesis, $G'$ has a vertex cover $A$ of size $m-1$. Then $A\cup \{a_1\}$ is a vertex cover of size $m$ in $G$. This completes the proof of Claim 2.

\medskip
\textbf{Claim 3.}~$\{a_{1},b_{m+1},c_{m+1}\},\{a_{m+1},b_{1},c_{m+1}\},\{a_{m+1},b_{m+1},c_{1}\}\notin E(G)$.
\medskip

Otherwise, without loss of generality, suppose that $\{a_{1},b_{m+1},c_{m+1}\}\in E(G)$.  By Claim 2, $\nu(G-a_1)=m$. Let $M$ be a matching of size $m$ in $G-a_1$ such that $V(M)\cap \{b_{m+1},c_{m+1}\}$ is minimal. We claim $V(M)\cap \{b_{m+1},c_{m+1}\}=\emptyset$. Otherwise, there exists $f\in V(M)$ such that $f\cap \{b_{m+1},c_{m+1}\}\neq \emptyset$. Without loss of generality, we may assume that $b_{m+1}\in f$. Since $m:=|M|$, there exists $b_j$ with $j<m+1$ such that $b_j\notin V(M)$. Recall that $G$ is partitely shifted. Hence we have $f'=(f\setminus \{b_{m+1}\})\cup \{b_j\}\in E(G)$. Thus $M'=(M\setminus\{f\})\cup \{f'\}$ is a matching of size $m$ in $G$ such that $|V(M')\cap \{b_{m+1},c_{m+1}\}|<|V(M)\cap \{b_{m+1},c_{m+1}\}|$, contradicting to the choice of $M$. Hence $M\cup \{a_1,b_{m+1},c_{m+1}\}$ is a matching of size $m+1$ in $G$, contradicting to that $\nu(G)=m$. This completes the proof of Claim 3.

Write $S=\{a_{m+1},b_{m+1},c_{m+1}\}$. Let $\mathcal{M}$ be a matching of size $m$ in $G$. Since $G$ is partitely shifted, we may assume that $V(M)=\bigcup_{i\in [m]}\{a_i,b_i,c_i\}$.

\medskip
\textbf{Claim 4.}~For any two edges $\{f_1,f_2\}\in {\mathcal{M}\choose 2}$, $L_S(f_1,f_2)\leq 4$.
\medskip

Otherwise, suppose there exists $\{f_1,f_2\}\in {\mathcal{M}\choose 2}$ such that $L_S(f_1,f_2)\geq 5$. By Lemma \ref{extension}, $L_S(f_1,f_2)$ has a matching size 3, say $\{g_1,g_2,g_3\}$, such that $g_1\in N_G(a_{m+1})$, $g_2\in N_G(a_{b+1})$ and $g_3\in N_G(c_{m+1})$. Then $(\mathcal{M}\setminus \{f_1,f_2\})\cup \{g_1\cup \{a_{m+1}\},g_2\cup \{b_{m+1}\},g_3\cup \{c_{m+1}\}\}$ is a matching of size $m+1$ in $G$, contradicting to that $\nu(G)=m$.

By Claims 2, 3 and 4, since $G$ is partitely shifted, one can see that
\begin{align*}
e(G)&\leq n^3-3m(n-m)^2-(n-m)^3-4(n-m){m\choose 2}\\
&=n^3-(n-m)(n^2+mn-2m)\\
&= 2mn+m^2n-2m^2,
\end{align*}
i.e.,
\begin{align}\label{up_e(G)}
e(G)\leq 2mn+m^2n-2m^2.
\end{align}
Other hand, since $b:=n-m\geq 2$ and $m\geq 2$, we have
\begin{align*}
&e(G)- 2mn+m^2n-2m^2\\
>&(m-1)n^2+2n-(2mn+m^2n-2m^2)\\
=&(m-1)(m^2+2mb+b^2)-(2mb+m^2(m+b))\\
=&m^2(b-1)+b^2(m-1)-4mb+2(m+b)\\
\geq& 2m(b-1)+2b(m-1)-4mb+2(m+b)= 0,
\end{align*}
contradicting to  (\ref{up_e(G)}). This completes the proof.\qed

\begin{theorem}[Daykin and H\"aggkvist, \cite{DH81}]\label{Daykin}
Let $n,k$ be two positive integers and let $H$ be a $k$-partite $k$-graph with $n$ vertices in each class.
If  $\delta_1(H)\geq 2n^{k-1}/3$, then $H$ has a perfect matching.
\end{theorem}

\begin{lemma}\label{matching-NPM}
Let $n\geq 162$ be an integer. Let $H$ be a 3-partite 3-graph with $n$ vertices in each class. If $e(H)>(n-2)n^2+2n$ and $\nu(H)=n-1$, then $\tau(H)=n-1$.
\end{lemma}
\pf By contradiction, suppose that $\tau(H)=n$. Since $H$ be a 3-partite 3-graph with $n$ vertices, then $G$ contains no isolated vertices, i.e.,  $\delta_1(H)\geq 1$. Next we will show that $H$ contains a perfect matching, which leads to a contradiction.

Let $B:=\{v\ |\ d_{H}(v)\leq 7n^2/9\}$. We claim $|B|\leq 9$.
Otherwise, suppose that $|B|\geq 10$. Without loss of generality, we may assume that $|B|=10$ since we may add some edges to $H$. Write $B_i=V_i\cap B$. Then we have
\begin{align*}
  e(H)&\leq n^3-|B|2n^2/9 +\sum_{\{i,j\}\in {[3]\choose 2}}n|B_i||B_j|-2|B_1||B_2||B_3|\\
  &\leq n^3-20n^2/9+33n \quad\mbox{(since $|B|\geq 10$)}\\
  &<(n-2)n^2+2n\quad \mbox{(since $n\geq 140$)}\\
  &<e(H),
\end{align*}
a contradiction.  Write $b:=|B|$ and $B:=\{v_1,\ldots,v_b\}$ such that $d_H(v_1)\leq \ldots\leq d_H(v_b)$.
Denote $B'=\{v_i\ |\ d_{H}(v_i)\leq 2bn\}$. We claim $|B'|\leq 2$, otherwise, we have 
\[
e(H) \leq n^3-3(n^2-2bn)+3n<(n-2)n^2+2n,
\]
a contradiction since $n\geq 56$ and $b\leq 9$.

We show that $H$ has a matching $M_1$ covering $B'$. If $|B'|\leq 1$, the result obviously holds since $\delta_1(H)\geq 1$. Now suppose that $B'=\{v_1,v_2\}$ such that $d_H(v_1)\leq d_H(v_2)$. If there exists an edge $f$ containing $\{v_1,v_2\}$, then the result is followed. So we may assume that $d_{H}(\{v_1,v_2\})=0$.  Since $e(H)\geq (n-2)n^2+2n+1$, then
\begin{align}\label{low-bound}
d_{H}(v_1)+d_H(v_2)\geq 2n+1,
\end{align}
which implies $d_H(v_2)\geq n+1$. Let $G_1:=N_H(v_1)$ and $G_2=N_{H}(v_2)$. Since $e(G_2)\geq n+1$ and $G$ is bipartite graph with at most $n$ vertices in each class, then $G_2$ contains two vertices $u_2,u_2'$  such that $\min \{d_{G_2}(u_2),d_{G_2}(u_2')\}\geq 2$. If $G_1-u_2$ has an edge $g_1$, then we may choose an edge $u_2\in g_2\in E(G_2\backslash V(g_1))$ such that   $\{g_1\cup\{v_1\},g_2\cup \{v_2\}\}$ is a matching covering $\{v_1,v_2\}$. So by symmetry, we may assume that $E(G_1-u_2)=\emptyset$ and $E(G_1-  u_2')=\emptyset$, which means that $E(G_1):=\{u_2u_2'\}$. If every edge of $G_2$ intersects $\{u_2,u_2'\}$, then $e(G_1)\leq 2n-1$, contradicting to (\ref{low-bound}). Thus $G_2-\{u_2,u_2'\}$ has an edge $g_2$. Then
$\{g_2\cup \{v_2\},\{v_1,u_2,u_2'\}\}$ is a desired matching covering $\{v_1,v_2\}$.

 Let $M_2$ be a matching of $H- V(M_1)$ such that every edge of $M_2$ intersects $B\backslash  B_1$ and $|V(M_2)\cap (B\backslash  B')|$ is maximum. We claim that $B\backslash B'\subseteq V(M_2)$.   Otherwise, suppose that there exists $v_i\in B\backslash (B'\cup V(M_2))$. Since $d_H(v_i)\geq 2bn$ and the number of edges containing $v_i$ and intersecting $V(M_1)\cup V(M_2)$ is less than $2bn$, then there exists an edge $f\in H-V(M_1\cup M_2)$ containing $v_i$. Thus $|V(M_2\cup \{f\})\cap (B\backslash  B')|>|V(M_2)\cap (B\backslash  B')|$, contradicting to the choice of $M_2$.

Let $H_1:=H\setminus V(M_1\cup M_2)$. Since $|M_1\cup M_2|\leq b$, we have
\[
\delta_1(H_1)\geq 7n^2/9-2bn\geq 2n^2/3,
\]
where the second inequality holds since $n\geq 162$.
By Theorem \ref{Daykin}, $H_1$ has a perfect matching $M_3$. Then $M_1\cup M_2\cup M_3$ is a perfect matching of $H$. This completes the proof. \qed


%

\section{Proof of Theorem~\ref{main-theorem}}

Since $n\geq 162$, by Lemma \ref{matching-small}, we may assume that $m\geq 13$.
By Lemma \ref{Stab-lem}, we may assume that $H$ is not  partitely shifted.  Let $H:=G_0,G_1,...,G_k$ be a sequence such that $G_{i+1}$ is obtained from $G_i$ by a shift operation and $G_k$ is  partitely shifted. By Lemma \ref{Stab-lem}, one can see that $\tau(G_k)=m$.

By contradiction, suppose $\tau(G)\geq m+1$. Then there exists $r\in [k]$ such that $\tau(G_r)\geq m+1$ and $\tau(G_{r+1})=m$. Suppose that $u,v\in V_i$ such that $G_{r+1}=S_{uv}(G_r)$. Let $W$ be a minimum vertex cover of $G_{r+1}$, then $|W|=m$. We claim that $u\in W$ and $v\notin W$. Otherwise, if either  $W\cap \{u,v\}=\emptyset$ or $\{u,v\}\subseteq W$ holds, then in both cases, $W$ is  a vertex cover of $G_r$, a contradiction. 
Moreover, one can see that $W\cup \{v\}$ is a vertex cover of $G_r$. Write $S:=W\cup \{v\}$. Note that $S$ is a minimum vertex cover of size $m+1$ in $G_r$. Without loss of generality, we may write that $H=G_r$ and $G:=G_{r+1}$. Then one can see that $G=S_{uv}(H)$, $e(H)=e(G)$ and $m=\nu(H)\geq \nu(G)$.

We discuss two cases.

\medskip
\textbf{Case 1.}~$m\leq \frac{3}{4}n$.
\medskip

Let $x\in {S\setminus \{u,v\}}$ such that $d_{H}(x)\leq d_{H}(y)$ for all $y\in {S\setminus \{u,v\}}$. Let $H'=H-({S\setminus\{u,v,x\}})$. Recall that $S$ is a vertex cover of $H$.  Since $e(H)\geq (m-1)n^2+3n-m$, then $d_{H'}(x)+d_{H'}(u)+d_{H'}(v)\geq n^2+3n-m$. Thus we have $max\{d_{H'}(x),d_{H'}(u),d_{H'}(v)\}\geq \frac{n^2}{3}.$
 Without loss of generality, suppose $d_{H'}(x)\geq \frac{n^2}{3}$. Note that $e(H'-x)\geq 3n-m\geq 2n+1$ and $\tau(H'-x)=2$. Let $\{x_1,x_2\}$ is a vertex cover of $H'-x$. One can see that $|N_{H'-x}(x_1)|+|N_{H'-x}(x_2)|>2n$.
 By Lemma \ref{matching-2size}, $\{N_{H'-x}(x_1),N_{H'-x}(x_2)\}$ admits a rainbow matching, which implies that
   $H'-x$ has a matching $M_0$ of size 2. Since the number of edges containing $x$ and intersecting $V(M_0)$ is at most $4n-4$ and $n>12$, there exists an edges $f\in H'-V(M_0)$ such that $x\in f$. Let $M_1:=M_0\cup \{f\}$.

Recall that $G=S_{uv}(H)$, $\tau(H)=m+1$ and $\tau(G)=m$. Thus for any $e\in N_{H}(v)$, $e\cap (W\setminus \{u\})\neq \emptyset$ or $e\notin N_H(u)$.  
For $i\in [3]$, write $B_i:=(S\backslash\{u,v,x\})\cap V_i$ and $A_i:=(V(M_1)\backslash \{u,x\})\cap V_i$. One can see that for $i\in [3]$, $|A_i|\geq 1$, $\sum_{i\in [3]}|A_i|=7$ and $\sum_{i\in [3]}|B_i|=m-2$.
If $d_{H}(x)\leq n^2-8n$,  since $W\setminus\{u,x\}$ is a vertex cover of $G-V(M_1)$, then we have
\begin{align*}
	e(H-V(M_1))&=e(G-V(M_1))\\
&\geq e(G)-d_{G}(u)-d_G(x)-e_G(W\backslash\{u,x\},V(M_1)\backslash\{u,x\})\\
&\geq (m-1)n^2+3n-m-n^2-(n^2-8n)\\
&-n\left((|A_1|+|A_2|)|B_3|-|A_1||A_2|+(|A_2|+|A_3|)|B_1|-|A_2||A_3|+(|A_1|+|A_3|)|B_2|-|A_1||A_3|\right)\\
&=(m-3)n^2+11n-m-n\left(\sum_{i\in[3]}|A_i|\sum_{i\in[3]}|B_i|-\sum_{i\in[3]}|A_i||B_i|-(|A_1||A_2|+|A_2|A_3|+|A_1||A_3|)\right)\\
&\geq (m-3)n^2+11n-m-n\left(7\sum_{i\in[3]}|B_i|-\sum_{i\in[3]}|B_i|-6\right)\\
	&= (m-3)n^2+11n-m-n(6(m-2)-6)\\
&=(m-3)(n-3)^2+11n-10m+27\\
	&> (m-3)(n-3)^2.
\end{align*}
Hence by Theorem \ref{Aharoni17}, $H-V(M_1)$ has a matching $M_2$ of size $m-2$. Then $M_1\cup M_2$ is a matching of size $m+1$ in $H$, a contradiction.

Now we may assume that $d_{H}(x)\geq n^2-8n$. Then for all $y\in {S\setminus\{u,v,x\}}$, we have $d_{H}(y)\geq n^2-8n$. Let $M_3$ be a maximal matching of $H-V(M_1)$ such that for  any $e\in M_3$, $|e\cap S|=1$. We claim that $|M_3|=m-2$. Otherwise, there exists $y\in S\backslash V(M_3\cup M_1)$. Note that the number of edges in  $H$ containing $y$ and intersecting with $(V(M_3)\cup S)$ is at most $n^2-(n-m)^2\leq 15n^2/16$. Since $n>128$, we have
\[
d_{H}(y)\geq n^2-8n>15n^2/16.
\]
Thus there exists an edge $f'$ containing $y$ in $H-V(M_1\cup M_3)$. Then $M_3':=M_3\cup\{f'\} $ is a matching of $H-V(M_1)$ such that for any $e\in M_3'$, $|e\cap S|=1$, contradicting to the choice of $M_3$. Now it follows that $M_1\cup M_3$ is a matching of size $m+1$ in $H$, contradicting to that $\nu(H)=m$.

\medskip
\textbf{Case 2.}~$m\geq \frac{3}{4}n$.
\medskip

\medskip
\textbf{Claim 2.}~$|W\cap (V_2\cup V_3)|\leq 1$.
\medskip

 Otherwise, suppose $|W\cap (V_2\cup V_3)|\geq 2$. Then
 \begin{align*}
 	e(H)=e(G) &\leq n^3-(n-m+2)(n-2)n\\
 	&=n(n^2-(n^2-mn+2(m-2)))\\
 	&=mn^2-2(m-2)n\\
 	&\le mn^2-2(\frac{3}{4}n-2)n\quad \mbox{(since $m\geq 3n/4$)}\\
 	&=(m-1)n^2-\frac{n^2}{2}+4n,\\
 	&< (m-1)n^2+3n-m,\\
 \end{align*}
 a contradiction, where the last inequality holds since $n>m\geq 1$. This completes the proof of Claim 2.

Next we discuss two subcases.

\medskip
\textbf{Subcase 2.1}~$|W\cap (V_2\cup V_3)|=0$.
\medskip

Then we have $\{u,v\}\subseteq V_1$ and $N_{H}(u)\cap N_{H}(v)=\emptyset$, which implies that  $d_{H}(u)+d_{H}(v)\leq n^2$. We choose $x\in S\setminus\{u,v\}$ such that $d_{H}(x)\leq d_{H}(y)$ for all $y\in S\setminus\{u,v\}$. Let $S_1=S\backslash\{u,v,x\}$ and $F=H-S_1$. Next we show that $\nu(F)\geq 3$. Note that $\{u,v,x\} $ is a vertex cover of $F$.   Since $e(F)\geq n^2+3n-m$, we have $d_{F}(u)+d_{F}(v)+d_{F}(x)\geq n^2+2n$. So we may choose that $x'\in \{u,v,x\}$ such that $d_{F}(x')\geq \frac{1}{3}n^2$. Note that $e(F-x')\geq 3n-m$ and $\tau(F-x')=2$.
Thus by Lemma \ref{matching-2size}, $F-x'$ has a matching of size 2. Let $\{e_1, e_2\}$ be a matching of $F-x'$. Since the number of edges containing $y$ and intersecting $e_1\cup e_2$ in $F$ is at most $4n-4<n^2/3$ ($n\geq 12$),  there exists an edge $e_3$ containing $x'$ in $F-V(\{e_1, e_2\})$. Thus $M_1=\{e_1, e_2, e_3\}$ is matching of size 3 in $H\setminus S_1$. Since $S_1$ is a vertex cover of $H-V(M_1)$, we have
\begin{align}\label{SC21-eq}
e_{H}(V(M_1)\backslash \{u,v,x\},S_1)\leq (m-2)(6n-9).
\end{align}
  So if $d_{H}(u)+d_{H}(v)+d_{H}(x)\leq 2n^2-4n+9$, then
  \begin{align*}
 	e(H-V(M_1))&\geq (m-1)n^2+3n-m-(d_{H}(u)+d_{H}(v)+d_{H}(x))-e_{H}(V(M_1)\backslash \{u,v,x\},S_1)\\
 &\geq (m-1)n^2+3n-m-(2n^2-4n+9)-(m-2)(6n-9)\quad \mbox{(by (\ref{SC21-eq}))}\\
 &> (m-3)(n-3)^2.
 \end{align*}
 By Theorem \ref{Aharoni17}, $H-V(M_1)$ has a matching of size $m-2$, say $M_2$. Then $M_1\cup M_2$ is a matching of size $m+1$ in $H$, contradicting to that $\nu(H)=m$.

   Next we may assume $d_{H}(u)+d_{H}(v)+d_{H}(x)> 2n^2-4n+10$. Since $d_{H}(u)+d_{H}(v)\leq n^2$, then $d_{H}(x)\geq n^2-4n+10$. Let $x\in e\in M_1$. Write $M_1':=M_1\backslash\{e\}$ and $H':=H-V(M_1')$.
Recall that for any $y\in S\setminus\{u,v\}$,
 \begin{align*}
 	d_{H-V(M_1')}(y)&\geq n^2-4n+10-(4n-4)\\
  	&=(n-2)^2-4(n-2)+2\\
  &>(n-6)(n-2).
 \end{align*}
We use $\{v_1,\ldots,v_{m-1}\}$ to denote $S-\{u,v\}$.
 If $m-2\leq n-6$, by Theorem \ref{Aharoni17-2}, $\{N_{H'}(v_i)\ |\ i\in [m-1]\}$ has a rainbow matching, which implies that $H'$ has a matching $M_2'$ of size $m-1$. It follows that $M_1'\cup M_2'$ is a matching of size $m+1$ in $H$,  a contradiction.

  By Lemma \ref{matching-NPM}, we may assume that $m\leq n-2$. Combining $m\geq n-3$,  we have $m\in \{n-2,n-3\}$. Let $t=n-m-1$. It follows that $t\in \{1,2\}$.  For $i\in [3]$, we choose $t$ vertices from  $V_i-V(M_1')$, saying $v_{i1},\ldots,v_{it}$ such that $\max\{d_{H'}(v_{ij})\ |\ j\in [t]\}\leq \min \{d_{H'}(y)\ |\ y\in V_i-(V(M_1')\cup \{v_{ij}\ |\ j\in [t]\})\}$. Write $T:=\cup_{i\in [3]}\{v_{ij}\ |\ j\in [t]\}$ and $H'':=H-(V(M_1')\cup T)$. If $\delta_1(H'')>\frac{2}{3}(n-t-2)^2$, then by Theorem \ref{Daykin}, $H''$ has a perfect matching of size $n-t-2$, which means that $H$ has a matching of size $m+1$, a contradiction. Hence we may assume $\delta_1(H'')\leq \frac{2}{3}(n-t-2)^2$.  We claim that $\delta_1(H'')> 4n$, otherwise, for some $i\in \{2,3\}$, there exist $t+1$ vertices $y\in V_i-V(M_1')$ such that $d_{H}(y)\leq 4n+(t+2)n$. So we have
  \begin{align*}
  (m-2)n^2+2n+2&\leq (m-2)n^2+3n-m\\
 &\leq  e(H-\{u,v\})\\
 &\leq (m-1)n^2-(t+1)(m-1)n+(4n+(t+2)n)(t+1)\\
 &\quad \mbox{(since $S\backslash\{u,v\}$) is a vertex cover of $H-\{u,v\}$}\\
 &\leq (m-1)n^2-2(m-1)n+14n,\quad \mbox{(since $m\geq 9$ and $1\leq t\leq 2$)}\\
 &\leq (m-1)n^2-2(n-4)n+14n,
  \end{align*}
  which implies that
   \begin{align*}
 n^2-20n+2\leq 0,
  \end{align*}
  contradicting to that $n\geq 20$.  Put $B:=\{y\ |\ d_{H''}(y)<2(n-t-2)^2/3+4n\}$. It follows that $|B|\leq 2$. Otherwise, if $|B|\geq 3$, let $B'\in {B\choose 3}$, we have $B'\subseteq V_2\cup V_3$ and so
   \begin{align*}
  &(m-2)n^2+2n+2\\
  &\leq  e(H-\{u,v\})\\
 &< (m-1)n^2-(t+|B'|)(m-1)n+(|B'|+t)(2(n-t-2)^2/3+4n)+(|B'|+t)(t+2)n\\
 &\leq (m-1)n^2-(t+3)(m-1)n+(t+3)(2(n-t-2)^2/3+4n)+(t+3)(t+2)n\\
 &\leq (m-1)n^2-4(m-1)n+4(2(n-t-2)^2/3+4n)+12n\\
 &\leq (m-1)n^2-4(n-4)n+8(n-3)^2/3+28n\\
 &= (m-2)n^2-n^2/3+28n+24,
  \end{align*}
  contradicting to that $n\geq 79$.

  Recall that  $\delta_1(H'')> 4n$ and $|B|\leq 2$. We may greedily find a matching $M_3$ of size $B$ in $H''$ such that $B\subseteq V(M_3)$. Moreover,
  \begin{align*}
    \delta_1(H''-V(M_3))> 2(n-t-2)^2/3+4n-2|B|n\geq 2(n-t-2)^2/3.
  \end{align*}
  Thus by Theorem \ref{Daykin},  $H''-V(M_3)$ has a perfect matching $M_4$ with size $n-t-2-|B|$, which means that $M_1'\cup M_3\cup M_4$ is a matching of size $m+1$ in $H$, a contradiction.

\medskip
\textbf{Subcase 2.2}~$|W\cap (V_2\cup V_3)|=1$.
\medskip

 Without loss of generality, suppose that $|W\cap V_2|=1$. Let $w\in {W\cap V_2}$ and  let $x\in {W\setminus \{u,w\}}$ such that $d_{G}(x)\leq d_{G}(y)$ for all $y\in W$.
Note that $x\in V_1$ and
\begin{align*}
 (m-1)n^2+3n-m&\leq  e(G)\\
&\leq (m-2)n^2+d_{G}(x)+d_{G\setminus\{W\setminus\ \{w\}\}}(w)\\
&\leq (m-2)n^2+d_{G}(x)+(n-m+1)n,
  \end{align*}
 which implies that
 $d_{G}(x)\geq mn+2n-m$.  By the definition of $G$, we know that $d_{H}(x)=d_{G}(x)\geq mn+2n-m$. Moreover,
 \begin{align*}
 d_{H-(W\backslash \{u\})}(u)+d_{H-(W\backslash \{u\})}(v)\geq e(H)-(m-1)n^2\geq 3n-m>2n.
 \end{align*}
Since $\{u,v\}$ is a vertex cover of $H-(W\backslash \{u\})$, by  Lemma \ref{matching-2size}, $H-(W\backslash \{u\})$ has a matching $\mathcal{M}_0$ of size two. The number of edges containing $x$ and intersecting $V(\mathcal{M}_0)$ in $H-(W\backslash \{u,x\})$ is at  most $4n-4$. Other hand, since $x\in V_1$ and $m\geq 13$, then the number of edges containing $x$ in $H-(W\backslash\{x\})$ is at least $mn+n-m>4n-4$. Thus there exists an edge $f\in H-((S\backslash \{x\})\cup V(\mathcal{M}_0))$ such that $\{x\}=f\cap S$. One can see that $\mathcal{M}_1:=\mathcal{M}_0\cup \{f\}$ is a matching of size three of $H-(S\backslash\{u,v,x\})$.

We first consider the case for $\{w\}=\{u\}$. 
Since $W$ is a vertex cover of $G$,
\begin{align*}
	e(H-V(\mathcal{M}_1)) &=e(G-V(\mathcal{M}_1))\\
&\geq e(G)-d_{G}(x)-d_{G}(u)+d_{G}(\{x,u\})-e_{G-\{u,x\}}(V(\mathcal{M}_1)\backslash (V_1\cup\{u\}),W\backslash\{u,x\})\\
&\geq
(m-1)n^2+3n-m-(2n^2-n)-(m-2)(5n-6)\\
&=(m-3)n^2-(5nm-5m-14n+12)\\
	&= (m-3)(n-3)^2+nm-4n-4m+15\\
	&> (m-3)(n-3)^2\quad \mbox{(since $m\geq 8$)}.
\end{align*}
Thus by Theorem \ref{Aharoni17}, $H-V(\mathcal{M}_1)$ has a matching of $\mathcal{M}_2$ size $m-2$. So $\mathcal{M}_1\cup \mathcal{M}_2$ is a matching of size $m+1$ in $H$, a contradiction.

Next we consider the case for $\{w\}\neq \{u\}$. It follows that $u,v\in V_1$. Since $W$ is a vertex cover,
 then we have
 \begin{align*}
d_{H-(W\backslash \{w\})}(w)&\geq d_{G-(W\backslash\{w\})}(w)\\
&\geq  e(G)-(m-1)n^2\\
&\geq 3n-m.
 \end{align*}
Other hand, 
\[
e_{H-(W\backslash\{w\})}(\{w\},V(\mathcal{M}_0))\leq  2(n-m)+(n-2)=3n-2m-2.
\]
Hence $H-((W\setminus \{w\})\cup V(\mathcal{M}_0))$ has an edge $f$ containing $w$  such that $|f\cap S|=1$. Write $\mathcal{M}_1':=\mathcal{M}_0\cup \{f\}$. Note that $W\setminus \{u,w\}$ is a vertex cover of $H-V(\mathcal{M}_1')$. Then we have
\begin{align*}
 e(H-V(\mathcal{M}_1'))&\geq e(H)-d_{H}(u)-d_{H}(v)-e_{H}(\{w\},V_1\setminus\{u,v\})-e_{H\setminus\{u,v,w\}}(W\setminus\{u,w\},V(\mathcal{M}_1')\setminus(V_1\cup \{w\}))\\
  &\geq e(H)-d_{H}(u)-d_{H}(v)-e_{H}(\{w\},V_1\backslash\{u,v\})-(m-2)(5n-6)\\
  &=e(G)-d_{G}(u)-e_{G}(\{w\},V_1\backslash\{u\})-(m-2)(5n-6)\\
  &\geq (m-1)n^2+3n-m-(2n^2-n)-(m-2)(5n-6)\\
  &=(m-3)(n-3)^2+(nm-4n-4m+15)\\
  &=(m-3)(n-3)^2+(n-4)(m-4)-1\\
  &>(m-3)(n-3)^2\quad \mbox{(since $n>m\geq 5$)}.
\end{align*}
  Thus by Theorem \ref{Aharoni17}, $H-V(\mathcal{M}_1')$ has a matching of $\mathcal{M}_2'$ size $m-2$. Then $\mathcal{M}_1'\cup \mathcal{M}_2'$ is a matching of size $m+1$ in $H$,  a contradiction.
 This completes the proof. \qed

\end{document}